\documentclass[a4paper]{llncs}

\usepackage{amssymb}
\usepackage{graphicx}
\usepackage{url}
\urldef{\mailsa}\path|{alfred.hofmann, ursula.barth, ingrid.haas, frank.holzwarth,|
\urldef{\mailsb}\path|anna.kramer, leonie.kunz, christine.reiss, nicole.sator,|
\urldef{\mailsc}\path|erika.siebert-cole, peter.strasser, lncs}@springer.com|

\newtheorem{thm}{Theorem} 
\newtheorem{defn}[thm]{Definition} 
\newtheorem{lem}[thm]{Lemma}

\usepackage{authblk}
\begin{document}

\mainmatter  

\title{CL-Shellable Posets with No EL-Shellings}

\titlerunning{CL-Shellable Posets with No EL-Shelling}

%
%
\author{Tiansi Li}

\pagestyle{plain}
\institute{Department of Mathematics and Statistics, Washington University in St. Louis}
\maketitle

%
%

\begin{abstract}

We construct an ungraded CL-shellable poset and a graded CL-shellable poset and show that neither is EL-shellable.
\end{abstract}

\section{Introduction}
Lexicographic shellability was introduced in the 1980s. Bj\"{o}rner defined EL-shellability in \cite{bjorner1980shellable}, and proved a conjecture of Stanley's that the existence of a labeling, satisfying certain conditions, of the edges of the Hasse diagram of a bounded, graded poset $P$ shows that $P$ is Cohen-Macaulay. Bj\"{o}rner and Wachs defined in \cite{bjorner1983lexicographically} the equivalent notions of CL-shellability and recursive atom ordering for a bounded, graded poset. They extended these notions to bounded, non-graded posets in \cite{bjorner1996shellable}. The theory of lexicographic shellability has been developed and applied by many authors, as one will see upon checking a list of papers that refer to \cite{bjorner1983lexicographically} and \cite{bjorner1996shellable}. While every EL-shellable poset is CL-shellable, the converse has remained open until now. In this paper, We present two examples where we show both ungraded and graded CL-shellable posets are not necessarily EL-shellable.

Before we give the proofs, we introduce some definitions and theorems. We recommend \cite{bjorner1983lexicographically}, \cite{bjorner1996shellable}, and \cite{wachs2006poset} for readers unfamiliar with lexicographic shellability.

\begin{defn} \cite{wachs2006poset}
For each face F of a simplicial complex $\Delta$, let $\langle F\rangle$ denote the subcomplex generated by $F$. A simplicial complex $\Delta$ is said to be shellable if its facets can be arranged in linear order $F_1$, $F_2$, $\dots$, $F_t$ such that the subcomplex $\left(\bigcup_{i=1}^{k-1}\left\langle F_{i}\right\rangle\right) \cap\left\langle F_{k}\right\rangle$ is pure and (dim$F_k$-1)-dimensional for all $k$= 2, $\dots$, $t$. Such an ordering of facets is called a shelling. A poset $P$ is shellable if its order complex $\Delta(P)$ is shellable.
\end{defn}

Shellability has been a very useful tool in topological combinatorics. Every shellable complex has the homotopy type of a wedge of spheres, and we can find the dimensions of the spheres given a shelling.

An edge-labeling of a bounded poset $P$ is a map from the edge set of the Hasse diagram of $P$ to a label poset $\Lambda$. We call a saturated chain in $P$ weakly increasing if the edge-labeling is weakly increasing in $\Lambda$ when we read the labels up along the chain.

\begin{defn} \cite[Definition 3.2.1]{wachs2006poset}
Let $P$ be a bounded poset. An edge-lexicographical labeling (EL-labeling, for short) of $P$ is an edge labeling such that in each closed interval $[x, y]$ of $P$ , there is a unique weakly increasing maximal chain, which lexicographically precedes all other maximal chains of $[x, y]$.
\end{defn}

We call $P$ EL-shellable if there is an EL-labeling of $P$.

\begin{thm} \cite[Theorem 3.2.2]{wachs2006poset}
Denote by $\overline{P}$ the poset obtained from $P$ by removing $\hat{0}$ and $\hat{1}$. Suppose $P$ is a bounded poset with an EL-labeling. Then the lexicographic order of the maximal chains of $P$ is a shelling of $\Delta(P)$. Moreover, the corresponding order of the maximal chains of $\overline{P}$" is a shelling of $\Delta(\overline{P})$.
\end{thm}

For any closed interval $[x, y]$, we call $[x, y]_r$ a rooted interval if $r$ is a maximal chain in $[\hat{0},x]$. A chain-edge labeling of a bounded poset is a map from the set of all pairs $(c, e)$ to the label poset $\Lambda$, where $c$ is a maximal chain of $P$ and $e$ is an edge in $c$, such that $(c, e)$ and $(c', e)$ get the same label if $c$ and $c'$ coincide from $\hat{0}$ to $e$. We obtain a label of a rooted edge $e_r = [x, y]_r$, where $x\lessdot y$, from the chain-edge label of $(c, e)$, where $c$ is a maximal chain containing $r$.

\begin{defn} \cite[Definition 3.3.1]{wachs2006poset}
Let $P$ be a bounded poset. A chain-lexicographic labeling (CL-labeling, for short) of $P$ is a chain-edge labeling such that in each closed rooted interval $[x, y]_r$ of $P$, there is a unique strictly increasing maximal chain, which lexicographically precedes all other maximal chains of $[x, y]_r$. A poset that admits a CL-labeling is said to be CL-shellable.
\end{defn}

Clearly, EL-shellability implies CL-shellability. And they both imply shellability.

\begin{thm} \cite[Proposition 2.3]{bjorner1983lexicographically}
EL-shellability $\Rightarrow$ CL-shellability $\Rightarrow$ Shellability.
\end{thm}

Recursive atom ordering is defined in \cite{bjorner1983lexicographically}, where it is shown that the existence of an Recursive atom ordering is equivalent to CL-shellability.

\begin{defn} \cite[Definition 4.2.1]{wachs2006poset}
A bounded poset $P$ is said to admit a recursive atom ordering if its length $l(P)$ is $1$ or if $l(P)>1$ and there is an ordering $a_1, a_2,\dots ,a_t$ of the atoms of $P$ that satisfies:
\begin{enumerate} 
    \item For all $j = 1,2,\dots,t$, the interval $[a_j, \hat{1}]$ admits a recursive atom ordering in which the atoms of $[a_j, \hat{1}]$ that belong to $[a_i, \hat{1}]$ for some $i<j$ come first.
    \item For all $i<j$, if $a_i, a_j < y$ then there is a $k<j$ and an atom $z$ of $[a_j, \hat{1}]$ such that $a_k<z\leq y$.
\end{enumerate}
A recursive coatom ordering is a recursive atom ordering of the dual poset $P^{*}$.
\end{defn}

\begin{thm} \cite[Theorem 4.2.2]{wachs2006poset}
A bounded poset $P$ is CL-shellable if and only if $P$ admits a recursive atom ordering.
\end{thm}

In the next two sections we present two examples CL-shellable posets that are not EL-shellable, one of which is ungraded and the other is graded.

\section{Ungraded Example}
We prove in this section that the Hasse diagram in figure 1 gives an example of an ungraded CL-shellable poset that does not admit any EL-shellings. The difference between CL-shellings and EL-shellings is that in a fixed interval $[x, y]$, CL-shellings allow different weakly increasing maximal chains of $[x, y]_r$ when we consider different roots $r$, whereas an in EL-shelling, there is a unique weakly increasing maximal chain in $[x, y]$ that does not depend on roots. In terms of recursive atom ordering, the difference is that the atom ordering above every element can be different depending on roots for a CL-shelling, while an EL-shelling induces a recursive atom ordering where the atom ordering above every element is independent of roots (see Lemma 8).

We claim that the poset as in figure 1 admits a recursive atom ordering where the atom order above each element, except at $y$, is independent of roots. And the recursive atom ordering in $[y, \hat{1}]_r$ relies on which root we choose. That is, the unique weakly increasing chain of $[y, \hat{1}]_r$ must be different when we consider the root through $a_1$ and the root through $a_6$. This implies that this poset cannot admit any EL-shelling.

\begin{figure}
    \centering
    \includegraphics[width=1\columnwidth]{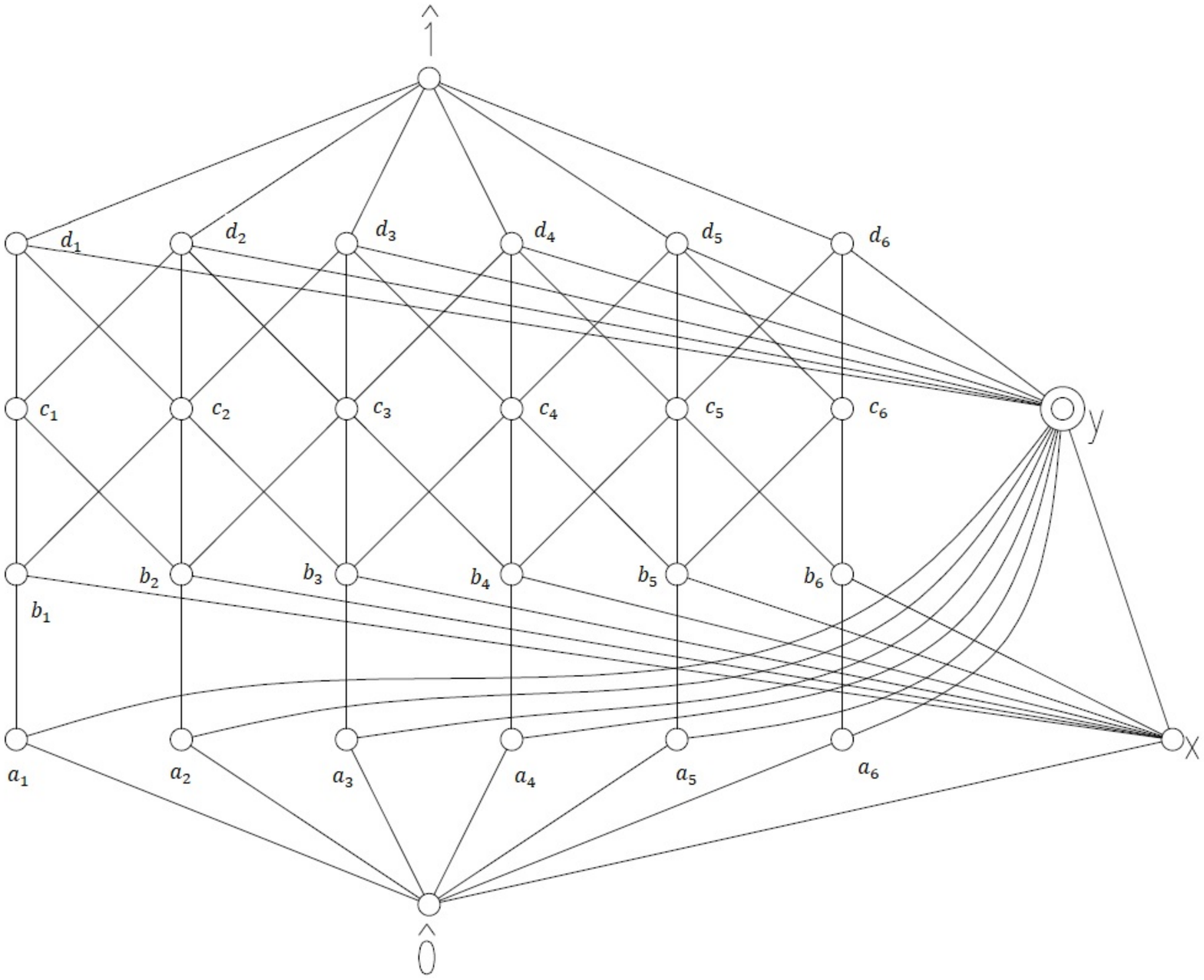}
    \caption{An ungraded CL-shellable poset with no EL-shelling}
    \label{fig:my_label}
\end{figure}

First we show that this poset admits a recursive atom ordering. Notice we only need to consider every element of height at most 2 except $y$, because for any element $z$ of height larger than 2 (or for $y$), every atom order of $[z, \hat{1}]$ induces a recursive atom ordering. We claim that listing from left to right, except at $\hat{0}$ where we list $x$ first and then every atom from left to right, gives a recursive atom ordering in every rooted interval.

For elements of height 2, notice from the diagram that $d_i$ is above all atoms of $[b_i, \hat{1}]$. Therefore, every ordering of the atoms of $[b_i, \hat{1}]$ satisfies the second condition in Definition 6. The recursive atom ordering in $[x, \hat{1}]$ follows similar arguments since for any two atoms $b_i$ and $b_j$ ($i<j$) of $[x, \hat{1}]$, a height 4 element above both atoms either covers a common atom of $[b_i, \hat{1}]$ and $[b_j, \hat{1}]$, or the leftmost atom of $[b_j, \hat{1}]$, in which case that atom of $[b_j, \hat{1}]$ covers some previous atoms above $x$. On the other hand, for elements $a_i$ of height 1, $[y, \hat{1}]$ has an atom that is above the other atoms of $[a_i, \hat{1}]$. Finally at $\hat{0}$, both atoms above $a_i$ cover $x$ for $i\in[6]$, and every element of height at least 2 is above $x$. Hence this poset admits a recursive atom ordering.

We state and prove a short lemma before finishing the proof.

\begin{lem}
An EL-shelling induces a recursive atom ordering in which the ordering of atoms above a given element does not depend on roots.
\end{lem}
Proof. First notice that any EL-shelling can be viewed as a CL-shelling where edge labels are independent of roots. By Theorem 7, a CL-shelling induces a recursive atom ordering in which the ordering of atoms above a given element with root $r$ is consistent with those edge labels with root $r$. That is, if for a fixed linear extension, the label of $[x, y]$ precedes the label of $[x, y']$, where $y$ and $y'$ both cover $x$, then $y$ precedes $y'$ in the atom ordering of $x$. Therefore if we start with an EL-shelling, the ordering of atoms above a given element in the induced recursive atom ordering does not depend on how one reached that elements from elements below it. \qed

Suppose we have an EL-shelling of the poset given in Figure 1. In the interval $[y, \hat{1}]$, there is an atom above $y$ that gives the unique increasing chain of that interval. It is independent of roots. Assume this atom is among $d_1$, $d_2$ and $d_3$. Consider the root of $y$ passing through $a_6$. Notice that $y$ is the second atom in its recursive atom ordering. So the atom above $y$ that comes first along this root must be among $d_4$, $d_5$ and $d_6$. This contradicts the assumption. If we now assume the atom that gives the unique increasing chain of $[y, \hat{1}]$ were among $d_4$, $d_5$ and $d_6$, the root of $y$ through $a_1$ gives a similar contradiction. Hence this poset cannot be EL-shellable.

\section{Graded Example}
We present in this section a graded poset that is CL-shellable but not EL-shellable.

The construction is based on a shellable but not extendably shellable complex exhibited by Hachimori. In \cite{hachimori2000combinatorics}, Hachimori constructed a shellable simplicial complex where the facet 134 comes last in every shelling of the complex (See Figure 2 below). By theorem 4.3 in \cite{bjorner1983lexicographically}, the dual of the face lattice of Hachimori's complex admits a recursive atom ordering. We will build a CL-shellable complex based on this poset, and we will show that the weakly increasing chain in $[134, \hat{1}]_r$ must be different for $r$ passing through $\hat{0}_a$ and $\hat{0}_d$. Therefore this poset cannot be EL-shellable.

\begin{figure}
    \centering
    \includegraphics[width=1\columnwidth]{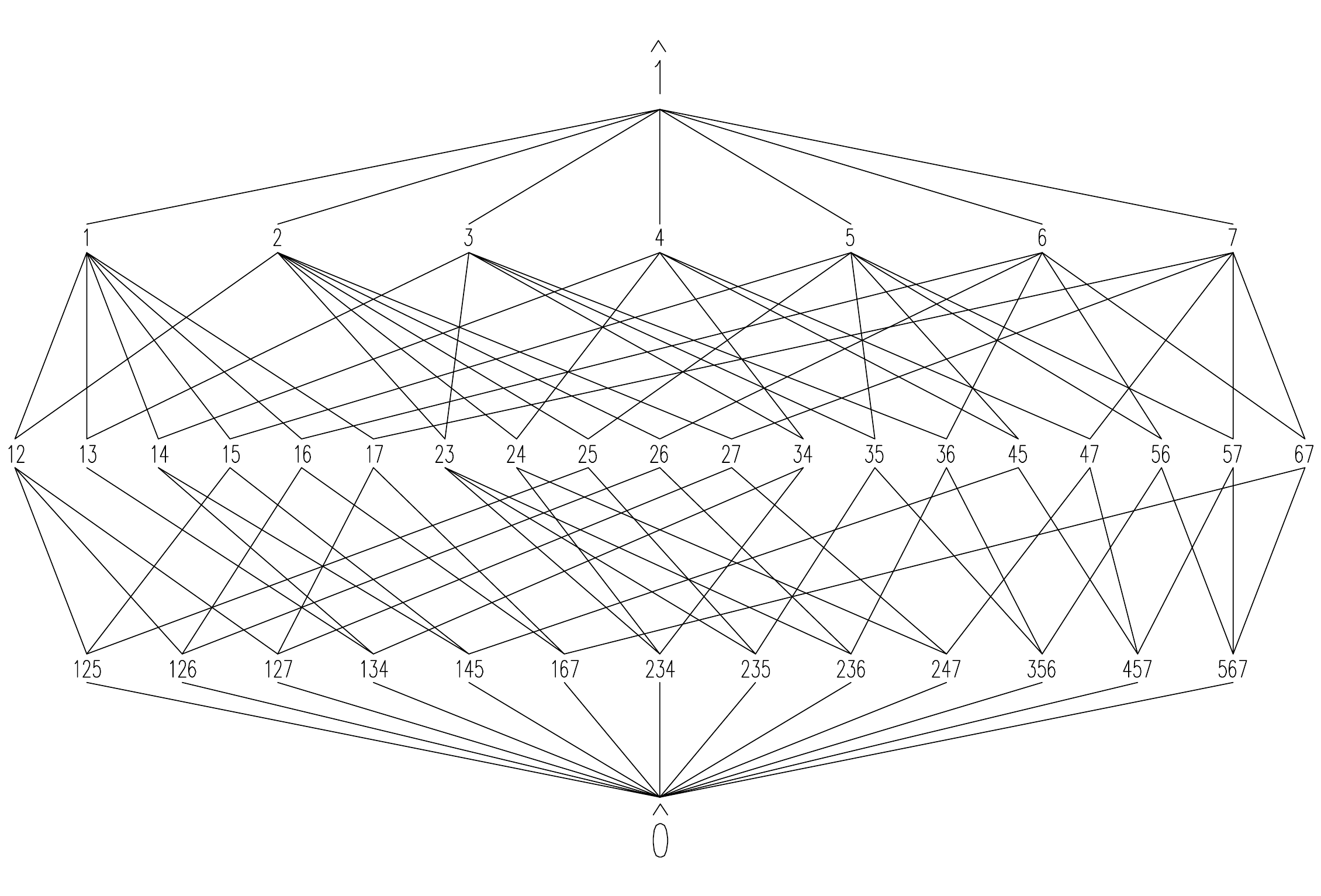}
    \caption{Dual of the face lattice of Hachimori's complex}
    \label{fig:my_label}
\end{figure}

Let us start with four copies of the dual of the face poset of Hachimori's complex. For convenience, we call them posets $A$, $B$, $C$, and $D$ and use $y_a$ or $y_b$ when refering to $y$ in $A$ or $B$ and so on. We build a new poset $P$ by first identifying all four copies of $134$ and $\hat{1}$ in $A$, $B$, $C$ and $D$, and then attaching a $\hat{0}$ to the bottom. Next we add a new element, call it $x$, which sits below all facet elements and above $\hat{0}$. So $x$ has rank 1. Finally we add an edge in the Hasse diagram between every pair of elements $y_i$ and $z_j$ if $y$ and $z$ represent non-empty faces in Hachimori's complex with $y$ a facet of $z$, and if $i$ and $j$ are consecutive in the lexicographic order. That is to say, for example, $y_b$ covers $z_a$, $z_b$ and $z_c$ whenever $y$ is a facet of $z$, but $y_b$ does not cover $z_d$, whereas $y_a$ covers $z_a$, $z_b$, but not $z_c$ or $z_d$. We claim that $P$ admits a recursive atom ordering, and $P$ admits no EL-shellings. Notice that $P$ is a graded poset of rank 5 with four copies of each element in the Hachimori's lattice except $134$ and $\hat{1}$.

Now we show that $P$ admits a recursive atom ordering. Fix a shelling order of the Hachimori's complex and a recursive coatom ordering of its face lattice induced by the shelling. Then each of A, B, C and D comes with a natural recursive atom ordering, from which we will build the recursive atom ordering of $P$. 

Suppose there is a recursive atom ordering of each $[\hat{0}_i, \hat{1}]$ and $[x, \hat{1}]$. Let the atom order for $[\hat{0}, \hat{1}]$ be $x$$\rightarrow$ $\hat{0}_a$$\rightarrow$ $\hat{0}_b$$\rightarrow$ $\hat{0}_c$$\rightarrow$ $\hat{0}_d$, since every atom of each $\hat{0}_i$ covers $x$, this induces a recursive atom ordering on $[\hat{0}, \hat{1}]$.

Suppose there is a recursive atom ordering in each of $[F_i, \hat{1}]$, where $F$ is a facet. Let the shelling order of Hachimori's complex be the atom order for $[\hat{0}_i. \hat{1}]$. If an element $y$ with index $i$ is above two atoms of $\hat{0}_i$, the recursive atom ordering of Hachimori's complex gives the existence of the element $z$, with index $i$, that satisfies the second condition in Definition 6. If the index of $y$ is not $i$, $z$ with index $i-1$ or $i+1$ satisfies the second condition in Definition 6.

For $[x, \hat{1}]$, we consider the following atom ordering:

\begin{center}
First facet with index $a$ in the shelling order

$\Downarrow$

First facet with index $b$ in the shelling order

$\Downarrow$

First facet with index $c$ in the shelling order

$\Downarrow$

First facet with index $d$ in the shelling order

$\Downarrow$

Second facet with index $a$ in the shelling order

$\Downarrow$

Second facet with index $b$ in the shelling order

$\Downarrow$

$\vdots$

$\Downarrow$

$134$
\end{center}

For any two atoms $y_i$, $z_j$ of $x$, the case where $y$=$z$ is obvious by construction. If $y$ is prior to $z$ in the shelling of Hachimori's complex, we have the following situations:

\begin{enumerate}
    \item If $i$ and $j$ are $a$ and $d$, an element above both $y_i$ and $z_j$ must be rank 4 with index $b$ or $c$. We can find an atom of $z_j$ with index $c$ below this element such that it covers some previous atoms of $x$.
    \item If $i$ and $j$ are $a$ and $c$, an element above both $y_i$ and $z_j$ is rank 4 with index $a$, $b$ or $c$, or is rank 3 with index $b$. In both cases we can find an atom of $z_j$ with index $b$ below this element (or it is the element in the rank 3 case), such that it covers some previous atoms of $x$.
    \item If $i$ and $j$ are $a$ and $b$, or $i=j$, we can find an appropriate atom of $z_j$ with index $i$ similarly.
    And all other cases are equivalent to one of the above.
\end{enumerate}

For rooted intervals $[F_i, \hat{1}]_i$ through $\hat{0}_i$, where $F$ is a facet in Hachimori's complex, the following atom order gives recursive atom ordering:

\begin{center}
Atoms with $i$-indices which cover facets prior to $F$

$\Downarrow$

Atoms with $j$-indices which cover facets prior to $F$, where $j$ and $i$ are consecutive letters

$\Downarrow$

Atoms with $i$-indices which do not cover facets prior to $F$

$\Downarrow$

All other atoms
\end{center}

For rooted intervals $[F_a, \hat{1}]_x$ through $x$, where $F$ is a facet in Hachimori's complex, the following atom order gives recursive atom ordering:

\begin{center}
Atoms with $a$-indices which cover facets prior to $F$

$\Downarrow$

Atoms with $b$-indices which cover facets prior to $F$

$\Downarrow$

Atoms with $a$-indices which do not cover facets prior to $F$

$\Downarrow$

Atoms with $b$-indices which do not cover facets prior to $F$
\end{center}

For rooted intervals $[F_i, \hat{1}]_x$ through $x$ where $i\neq a$ and $F$ is a facet in Hachimori's complex, the following atom order gives recursive atom ordering:

\begin{center}
Atoms with $i-1$-indices which cover facets prior to $F$

$\Downarrow$

Atoms with $i$-indices which cover facets prior to $F$

$\Downarrow$

Atoms with $i+1$-indices (if they exist) which cover facets prior to $F$

$\Downarrow$

Atoms with $i+1$-indices (if they exist) which do not cover facets prior to $F$

$\Downarrow$

All other atoms
\end{center}

Notice that we can break ties arbitrarily in the process described above because Hachimori's complex is a simplicial complex.

As for rooted intervals $[134, \hat{1}]_i$, we can still follow those steps except $i$ stands for the index of $\hat{0}_i$ (root).

For the rooted interval $[134, \hat{1}]_x$ through $x$, we order the atoms as: $14_a \rightarrow 34_a \rightarrow 14_b \rightarrow 34_b \rightarrow \dots \rightarrow 34_d \rightarrow 13_a \rightarrow \dots \rightarrow 13_d$

For the remaining elements $z$, the length of $[z, \hat{1}]$ is at most two, hence every atom order induces a recursive atom ordering. We have shown that $P$ admits a recursive atom ordering.

Consider the interval $[134, \hat{1}]$. It is a rank 3 interval with 12 atoms. Suppose $P$ admits an EL-shelling, where the unique increasing chain in this interval goes through one of $13_a$, $14_a$, $34_a$ or $13_b$, $14_b$, $34_b$. Consider the root of $134$ through $\hat{0}_d$. None of the six atoms with $a$ or $b$-indices cover any atoms of $\hat{0}_d$ other than $134$, whereas each of the six atoms with $c$ or $d$-indices cover some atoms of $\hat{0}_d$ other than $134$. Since $134$ is last in any recursive atom ordering, every atom of $134$ with $c$ or $d$-indices must be prior to every atom with $a$ or $b$-indices, and we have a contradiction. Similarly we can get a contradiction by assuming the EL-shelling implies one of the atoms with $d$-indices being increasing chain and take the root through $\hat{0}_a$. Hence $P$ cannot be EL-shellable.

\section*{Remark}
Notice that both examples in this paper require particular models to begin with. That is, a CL-shellable poset in which there exists two atoms $a$ and $b$ such that $a$ is prior to $b$ in every recursive atom ordering. We used the poset consisting of a 3-chain and a 2-chain (a pentagon in the Hasse diagram) in the ungraded example and the dual of the face lattice of Hachimori's complex in the graded example, both of which have an atom that must come last in any recursive atom ordering. It remains open whether one can construct a CL-shellable but not EL-shellable poset such that for any element $e$ and any two atoms $a$, $b$ of $e$, there exists two recursive atom orderings on $[e, \hat{1}]$ where $a$ is prior to $b$ in one and $b$ is prior to $a$ in the other.

\section*{Acknowledgement}
The author thanks Quancheng Mu and Moya Xiong for CAD plotting, and Russ Woodroofe for providing Hachimori's complex as a foundation the graded example.

\bibliographystyle{plain}
\bibliography{bibliography.bib}

\end{document}